\documentclass[12pt]{amsart}

\usepackage{latexsym,amsmath,amssymb,amsthm,amsfonts}
\usepackage{hyperref}

\usepackage[lite,alphabetic]{amsrefs}

\newcommand{\R}{{\mathbb{R}}}

\newcommand{\V}{{\mathcal{V}}}
\newcommand{\W}{{\mathcal{W}}}

\newcommand{\G}{{\widetilde G}}
\newcommand{\M}{{\widetilde M}}

\renewcommand{\phi}{\varphi}

\def\be{\begin{equation}}
\def\ee{\end{equation}}

\newtheorem{theorem}{Theorem}%[section]

\theoremstyle{remark}

\numberwithin{equation}{section}

\title{There is no strongly regular graph with parameters $(460, 153,32,60)$ }

\author{A. V. Bondarenko}
\address{Department of Mathematical Sciences, Norwegian University of Science and Technology, NO-
7491 Trondheim, Norway \newline and\newline
Department of Mathematical Analysis, National Taras Shevchenko
University, str. Vo\-lo\-dy\-myr\-ska, 64, Kyiv, 01601, Ukraine}
\email{andriybond@gmail.com}

\author{A. Mellit}
\address{International School for Advanced Studies (SISSA), Via Bonomea 265, 34136, Trieste,
Italy
\newline and\newline
International Centre for  Theoretical Physics, Trieste, Italy}
\email{mellit@gmail.com}

\author{A. Prymak}
\address{Department of Mathematics, University of Manitoba, Winnipeg, MB, R3T2N2, Canada}
\email{prymak@gmail.com}

\author{D. Radchenko}
\address{International Centre for  Theoretical Physics, Trieste, Italy}
\email{danradchenko@gmail.com}

\author{M. Viazovska}
\address{Humboldt University of Berlin, Berlin Mathematical School,  Rudower Chaussee 25, 12489 Berlin, Germany}
\email{viazovska@gmail.com}

\thanks{}

\keywords{Strongly regular graph, Euclidean representation}

\subjclass[2010]{Primary 05C25. Secondary 05C50, 52C99, 41A55.}

\begin{document}

\begin{abstract}

We prove that there is no strongly regular graph (SRG) with parameters $(460,153,32,60)$. The proof is based on a recent lower bound on the number of 4-cliques in a SRG and some applications of Euclidean representation of SRGs.
\end{abstract}

\maketitle

\section{Introduction}

%[[More historical info]]

A finite, undirected, simple graph $G=(V,E)$ with vertices $V$ and edges $E$ is called \emph{strongly regular} with parameters $(v,k,\lambda,\mu)$ if $G$ is $k$-regular on $v$ vertices, and, in addition, any two adjacent vertices of $G$ have exactly~$\lambda$ common neighbors, while any two non-adjacent vertices of $G$ have exactly~$\mu$ common neighbors.

The parameters~$(v,k,\lambda,\mu)$ of a $SRG$ must satisfy certain known 
conditions (see~\cite{Br-vL-82}), but in general it is an open 
question to determine parameters $(v,k,\lambda,\mu)$ for which strongly regular 
graphs (SRGs) exist, and, in case when they do exist, to classify such graphs. 
A list of known results for $v\le 1300$ is maintained at~\cite{Br-www}.

As an application of a recently established lower bound on the number of $4$-cliques in a SRG (see~\cite{BPR} and~\cite{BPR-web}) and of Euclidean representation of SRGs, we obtain the following non-existence result in a very special case.
\begin{theorem}\label{main}
There is no strongly regular graph with parameters $(460,153,32,60)$.
\end{theorem}

Some general background on SRGs can be found in~\cite[Chapter~9]{Br-Ha} 
and~\cite{Cam}, while~\cite[Chapter~8]{Br-Ha} and~\cite{Br-vL-82} contain 
details on Euclidean representation of SRGs. The argument with the Gram matrix 
used here has been extensively applied in~\cite{BPR}.

%Let us suggest the book~\cite{Br-Ha} as a reference for background needed in 
%this note. In particular, \cite[Chapter~9]{Br-Ha} and~\cite{Cam} overview 
%various properties of SRGs, while~\cite[Chapter~8]{Br-Ha} and~\cite{Br-vL-82} 
%contain some details on Euclidean representation of SRGs. The argument with %the Gram matrix used here has been extensively applied in~\cite{BPR}.

\section{Proof of Theorem 1}
Assume that a SRG $G=(V,E)$ with parameters $(v,k,\lambda,\mu)=(460,153,32,60)$ 
exists. Then adjacency matrix $A$ of $G$ satisfies the equations $AJ=153J$ and 
$A^2+28A-93I=60J$, where $I$ is the identity matrix and $J$ is the matrix 
having all entries equal to $1$. Consequently, $A$ has the following spectrum: 
$153^1\,3^{414}\,(-31)^{45}$. The Euclidean representation of $G$ defines a 
mapping $V\ni u\mapsto x_u \in\R^{45}$ such that all $x_u$ are unit vectors
and the dot products between these 
vectors depend only on the adjacency between the corresponding vertices of $G$. 
More precisely, for two different vertices $u,w\in V$, we have
	\[\langle x_u, x_w \rangle =
	\begin{cases}
	p, & \text{if }u\text{ is adjacent to }w,\\
	q, & \text{if }u\text{ is not adjacent to }w,\\
	\end{cases}
	\qquad\text{where}\quad
	\begin{array}{l}
	p=\frac{-31}{153},\\
	q%=\frac{15}{153}
	=\frac{5}{51},
	\end{array}\]
and $\langle x,y\rangle$ is the Euclidean dot product in $\R^{45}$. 
%The standard way to construct $x_u$ is to define them as rows (or columns) of a suitable linear combination of $A$, $I$, and $J$.

The first step is to show that $G$ has at least $228111$ complete 
subgraphs of size~$4$. This follows from~\cite{BPR} (see also the
bounds on the number of $4$-cliques in~\cite{BPR-web}). Let us give
a brief outline of the proof. For each edge $e=\{u,w\}\in E$ 
we consider the unit vector $y_e = \frac{x_u+x_w}{\|x_u+x_w\|}$. A simple 
calculation shows that the \textit{distribution} of dot products 
between $\{x_u\}_{u\in V}$ and $\{y_e\}_{e\in E}$ depends only on the 
SRG parameters and the number of $4$-cliques. Since the Gegenbauer 
polynomials $C_t^{(d-2)/2}(x)$ are positive definite on $S^{d-1}$, by applying
them to the Gram matrix of $\{x_u\}_{u\in V}\cup \{y_e\}_{e\in E}$ we get 
a positive definite matrix (here $d=45$, $t=4$). 
To get an inequality for the number of $4$-cliques we simply compute the value 
of the corresponding quadratic form on the vector that takes value $1$ 
on $x_u$'s and $a$ on $y_e$'s and optimize the parameter $a$.

Next, for any two adjacent vertices $u, w\in V$ (i.e., $\{u, w\}\in E$), 
let $V_{u,w}$ be the set of vertices $t\in V$ adjacent to both $u$ and $w$. Note that any pair of adjacent vertices in $V_{u,w}$ forms a $4$-clique together with $u$ and $w$. Now choose $u, w$ so that the number of edges in the subgraph of $G$ induced by $V_{u,w}$ is largest possible (among all possible edges $\{u, w\}\in E$). Let $\V:=V_{u,w}$, $\G$ be the subgraph of $G$ induced by $\V$, and let $m$ be the number of edges in $\G$. Since $G$ has $35190$ edges, we get the following inequality on $m$ from the lower bound on the number of $4$-cliques: $m\ge \frac{6\cdot 228111}{35190}$, so $m\ge 39$. For an upper bound on $m$, we will make use of the above Euclidean representation. Define $X_1:=\sum_{t\in \V}x_t$ and $X_2:=x_u+x_w$. The Gram matrix $M:=(\langle X_i, X_j\rangle)_{i,j=1}^2$ is positive semi-definite, therefore, $\det M \ge 0$. 
%We can compute the entries of $M$ using $m$ and parameters of our graph, namely
Explicitly, in terms of $m$ and graph parameters we have
	\[
	M=\begin{pmatrix}
	\lambda + 2mp + (\lambda^2-\lambda-2m)q & 2\lambda p \\
	2 \lambda p & 2+2p
	\end{pmatrix}=\frac{1}{153^2}
	\begin{pmatrix}
	19776-92m & -1984 \\
	-1984 & 244
	\end{pmatrix}.
	\]
The inequality $\det M\ge 0$ leads to $m \le \frac{2416}{61}$, therefore $m\le 39$.

Thus, $m=39$, and the graph $\G$ on $32$ vertices has $39$ edges. Let $\W$ be a set of $14$ vertices of $\G$ which have the largest degrees (in $\G$). We claim that the sum $\alpha$ of the degrees of vertices of $\W$ in $\G$ is at least $42$. Indeed, if each such degree is at least $3$, then we are clearly done. Otherwise, the sum of the degrees of vertices not in $\W$ is at most $(32-14)2=36$, which means that $\alpha\ge 2\cdot 39-36=42$. Denote by $\beta$ the number of edges in the subgraph induced by~$\W$. Then we have $\alpha-2\beta$ edges between $\W$ and $\V\setminus \W$, and $39+\beta-\alpha$ edges in $\V\setminus \W$. We take $Y_1:=\sum_{t\in \V\setminus \W}x_t$, $Y_2:=\sum_{t\in \W}x_t$, and $Y_3:=x_u+x_w$ and apply previous considerations. For the Gram matrix $\M:=(\langle Y_i, Y_j\rangle)_{i,j=1}^3$ we clearly have $\det \M\ge0$. On the other hand, we compute
% which can be expressed in terms of $\alpha$ and $\beta$. More specifically,
\begin{align*}
\langle Y_1,Y_1 \rangle &= 18+2(39+\beta-\alpha)p+(18\cdot 17-2(39+\beta-\alpha))q, \\
\langle Y_2,Y_2 \rangle &= 14+2\beta p + (14\cdot13-2\beta)q, \\
\langle Y_1,Y_2 \rangle &= (\alpha-2\beta)p+(18\cdot14-(\alpha-2\beta))q, \\
\langle Y_1,Y_3 \rangle &= 18\cdot2p, \quad
\langle Y_2,Y_3 \rangle = 14\cdot2p, \quad
\langle Y_3,Y_3 \rangle = 2+2p,
\end{align*}
and therefore
\[
\det\M=\left(-\frac{516304}{3581577}\alpha^2 + \frac{35785792}{3581577}\alpha\right)  - \left( \frac{1252672}{3581577}\beta +
\frac{198599296}{1193859}\right) =: \phi(\alpha) - \psi(\beta).
\]
The quadratic function $\phi$ is decreasing for $\alpha\ge \frac{35785792}{2\cdot 516304}=\frac{2114}{61}$, in particular for $\alpha\ge 42$. The linear function $\psi$ is clearly increasing. Since $39+\beta-\alpha\ge0$, we have $\beta\ge3$. Now, since
	\[0\le \det\M=\phi(\alpha) - \psi(\beta)\le \phi(42) - \psi(3)=\frac{-270848}{132651}<0,\]
we get a contradiction and hence Theorem~1 is proved.

\section{Conclusion}
Let us remark that the exact same reasoning from the proof above can be applied
to some other strongly regular graphs. For instance, with some trivial changes 
we obtain non-existence of strongly regular graphs with parameters 
$(5929,1482,275,402)$ and $(6205,858,47,130)$, for which the number 
of~$4$-cliques is bounded from below by~$4805$ and~$113$ respectively. The key 
property that these three graphs have in common is that they have a very 
small (but strictly positive) value of the Krein parameter $q_{22}^2$. The 
proof also goes through for some strongly regular graphs that 
satisfy $q_{22}^2 = 0$, or equivalently
	\[(s+1)(k+s+rs) = (k+s)(r+1)^2.\] 
In this case the above reasoning shows that all $\lambda$-subgraphs must 
be regular. The smallest set of parameters that can be ruled out in this way 
is $(2950,891,204,297)$.
Alternatively, the non-existence in this case can be shown by noting that the 
first subconstituent must be strongly regular, but there exist no strongly 
regular graphs on~$891$ vertices of degree~$204$ (since there are no 
feasible parameters with $v=891$ and $k=204$).

\end{document}